\documentclass{article}

\title{On singular univariate specializations of bivariate hypergeometric functions}

\author{Raimundas Vid\=unas\\
 \em Kobe University}

\newtheorem{theorem}{Theorem}[section]

\newtheorem{lemma}[theorem]{Lemma}

\usepackage{latexsym}

\newcommand{\app}[4]{F_{\!#1}\!
  \left(\left.{#2 \atop #3}\right| #4 \right) }

\newcommand{\hpg}[5]{{}_{#1}\mbox{\rm F}_{\!#2}\!
  \left(\left.{#3 \atop #4}\right| #5 \right) }

\newcommand{\hpgo}[2]{{}_{#1}\mbox{\rm F}_{\!#2}}

\newcommand{\proof}{{\bf Proof. }}
\newcommand{\qed}{\hfill $\Box$}
\newcommand{\equal}{&\!\!\!=\!\!\!&}

\begin{document}

\maketitle

\begin{abstract} 
It is tempting to evaluate $F_2(x,1)$ 
and similar univariate specializations of Appell's functions 
by evaluating the apparent power series at $x=0$ straight away
using the Gauss formula for ${}_2F_1(1)$.
But this kind of naive evaluation can lead to errors
as the ${}_2F_1(1)$ coefficients might eventually diverge;
then the actual power series at $x=0$ might involve branching terms.
This paper demonstrates these complications by concrete examples.
\end{abstract}

\section{Introduction}
\label{sec:symsq}

The paper \cite{AppelGhpg} investigated the question of when univariate specializations
of Appell's hypergeometric functions satisfy the same ordinary Fuchsian equation as familiar univariate hypergeometric functions. Several such coincidences of ordinary Fuchsian equations were found,
and then attempts were made to relate the Appell's and usual univariate hypergeometric solutions directly. In particular, the following wrong formulas were claimed in a neighborhood of $x=0$:
\begin{eqnarray} \label{we:ap2y1}
\app2{a;\;b_1,b_2}{c_1,c_2}{x,\,1}\equal
\!\frac{\Gamma(c_2)\Gamma(c_2-a-b_2)}{\Gamma(c_2-a)\Gamma(c_2-b_2)}
\hpg32{a,\,b_1,\,a-c_2+1}{\!c_1,a+b_2-c_2+1}{x},\\
\label{we:ap2xy1} \app2{a;\;b_1,b_2}{c_1,c_2}{x,\,1-x}\equal
\frac{\Gamma(c_2)\Gamma(c_2-a-b_2)}{\Gamma(c_2-a)\Gamma(c_2-b_2)}\,
(1-x)^{-a}\times \nonumber\\
&& \hpg32{a,\,c_1-b_1,\,a-c_2+1}{c_1,\,a+b_2-c_2+1}{\frac{x}{x-1}},\\
\label{we:ap4f3} \app4{a;\,b}{\!c+\frac12,\frac12}{x^2,(1\!-\!x)^2\!}\! \equal
\frac{\Gamma(\frac12)\Gamma(\frac12-a-b)}{\Gamma(\frac12-a)\Gamma(\frac12-b)}\,
\hpg32{2a,\,2b,\,c}{a+b+\frac12,2c}{x},\\
\label{we:ap4f2} \app4{\!2c\!-\!\frac12;3c\!-\!1}{c+\frac12,c+\frac12}{x^2,(1\!-\!x)^2\!}\!\! \equal
\frac{\Gamma(c+\frac12)\Gamma(2-4c)}{\Gamma(1-c)\Gamma(\frac32-2c)}\,
\hpg21{c,\,3c-1}{2c}{x}^2,
\end{eqnarray}
under the following conditions, respectively: $\mbox{Re}(c_2-a-b_2)>0$,
$\mbox{Re}(c_2-a-b_2)>0$, $\mbox{Re}(a+b)<\frac12$, $\mbox{Re}(c)<\frac12$.
On the left-hand side of these formulas, we have univariate specializations of the
bivariate Appell's hypergeometric functions $F_2(u,v)$ or $F_4(u,v)$. 
They are restricted to the singular curves $v=1$, $u+v=1$ or $\sqrt{u}+\sqrt{v}=1$
of the respective holonomic systems of partial differential equations for Appell's $F_2$ or $F_4$.
As well known \cite[Table 1]{Srivastava85}, Appell's $F_2$ or $F_4$ series converge
in the regions $|u|+|v|=1$ or $\sqrt{|u|}+\sqrt{|v|}=1$, respectively. 
Evaluation at the points  $x\in[0,1)$ of the Appell functions
in (\ref{we:ap2xy1})--(\ref{we:ap4f2}) is supposed to give the limiting values from 
inside the relevant convergence regions. The evaluation of the $F_2$ function in (\ref{we:ap2y1})
is supposed to apply to a one-dimensional restriction on a two-dimensional branch 
originating from the corner $(u,v)=(0,1)$ of the convergence region. The ``formulas"
(\ref{we:ap2y1})--(\ref{we:ap4f2}) are discussed in \cite{AppelGhpg} in two paragraphs:
in the paragraph containing formula (32) there, 
and the last paragraph of Section 4.

The rationale for each of the ``formulas" (\ref{we:ap2y1})--(\ref{we:ap4f2}) is that both sides satisfy the same Fuchsian ordinary equation (of order 3), both sides are proper power series at $x=0$ (and the linear space of solutions of the Fuchsian equation with integer local exponents at $x=0$ is one-dimensional in general), and that the $\Gamma$-factor is correct due to evaluation by
Gauss' formula \cite[Theorem 2.2.2]{specfaar} at $x=0$. However, the impression that the
restricted Appell's functions behave as proper power series at $x=0$ is wrong. In particular,
we can write the $F_2$ series in (\ref{we:ap2y1}) as follows:
\begin{equation} \label{eq:app2x0}
\app2{a;\;b_1,b_2}{c_1,c_2}{x,\,1}=\sum_{k=0}^{\infty} \frac{(a)_k\,(b_1)_k}{(c_1)_k\,k!}\,
\hpg21{a+k,\,b_2}{c_2}{\,1}x^k.
\end{equation}
But the convergence condition for the $\hpgo21(1)$ series is $\mbox{Re}(c_2-a-b_2-k)>0$,
so the coefficients in the power series are undefined for large enough $k$.
Quite similarly, the $k$th derivative of the $F_2$ function in (\ref{we:ap2y1}) at $x=0$ 
would evaluate to a linear combination of the values
\begin{equation}
\hpg21{a+k,\,b_2+j}{c_2+j}{\,1\,}, \qquad j=0,1,\ldots,k.
\end{equation}
These $\hpgo21(1)$ values are undefined for larger enough $k$ as well.
The specialized Appell's series could have branching behavior at $x=0$, 
and ``formulas" (\ref{we:ap2y1})--(\ref{we:ap2xy1}) must be generally wrong.
The main contribution of this paper is to demonstrate this branching behavior
on a concrete example.

The specialized $F_4$ series in (\ref{we:ap4f3})--(\ref{we:ap4f2}) has the same problem:
the terms in the ``obvious" power series expansion at $x=0$ are undefined for high enough degrees.
On the other hand, the following two formulas \cite[(40) and (76)]{AppelGhpg} are correct
if $\mbox{Re}(c-a_2-b_2)>0$ or $\mbox{Re}(c-a-b_2)>0$, respectively:
\begin{eqnarray}
\app3{a_1,a_2;\;b_1,b_2}{c}{x,\,1} \equal
\frac{\Gamma(c)\Gamma(c-a_2-b_2)}{\Gamma(c-a_2)\Gamma(c-b_2)}\,
\hpg32{a_1,b_1,c-a_2-b_2}{c-a_2,c-b_2}{x},\\
\app1{a;\;b_1,b_2}{c}{x,\,1}\equal
\frac{\Gamma(c)\Gamma(c-a-b_2)}{\Gamma(c-a)\Gamma(c-b_2)}\,
\hpg21{a,\,b_1}{c-b_2}{x}.
\end{eqnarray}
Here we can expand the $F_3$ or $F_1$ series as
\begin{eqnarray} \label{eq:f31}
\app3{a_1,a_2;\;b_1,b_2}{c}{x,\,1}\equal \sum_{k=0}^{\infty} \frac{(a_1)_k\,(b_1)_k}{(c)_k\,k!}\,
\hpg21{a_2,\,b_2}{c+k}{\,1}x^k,\\ \label{eq:f11}
\app1{a;\;b_1,b_2}{c}{x,\,1}\equal \sum_{k=0}^{\infty} \frac{(a)_k\,(b_1)_k}{(c)_k\,k!}\,
\hpg21{a+k,\,b_2}{c+k}{\,1}x^k.
\end{eqnarray}
The conditions for absolute convergence of the $\hpgo21(1)$ series are $\mbox{Re}(c-a_2-b_2+k)>0$
or $\mbox{Re}(c-a-b_2)>0$, respectively.  Therefore we can apply Gauss' summation formula
and take term-wise limits if only $\mbox{Re}(c-a_2-b_2)>0$ or $\mbox{Re}(c-a-b_2)>0$, respectively.  
For $x\in(-1,1)$, the specialized $F_3$ or $F_1$ series take limit values on the edge
of the convergence region $\max(|u|,|v|)<1$ for the $F_3(u,v)$ or $F_1(u,v)$ series.
In Section \ref{sec:appell} we summarize wrong and correct formulas given in \cite{AppelGhpg}.

\section{An explicit example}

Here we present our main counter-example to formulas (\ref{we:ap2y1})--(\ref{we:ap2xy1}).
The following Appell's $F_2$ function is under consideration:
\begin{lemma}
For $a\not\in\{1,2\}$ and $|x|+|y|<1$, 
\begin{equation} \label{eq:specf2}
\app2{a;1,1}{2,2}{x,y}=\frac{1-(1-x)^{2-a}-(1-y)^{2-a}+(1-x-y)^{2-a}}{(1-a)(2-a)\,x\,y}
\end{equation}
\end{lemma}
\proof We start with the following univariate identity:
\begin{equation}
\hpg21{a,\,1\,}{2}{x}=\sum_{i=0}^{\infty} \frac{(a)_i\,x^i}{(i+1)!}=\frac{1-(1-x)^{1-a}}{(1-a)\,x}.
\end{equation}
To show this formula, recognize the series of $1+(a-1)\,x\,\hpgo21(x)$ as the Taylor series for
$(1-x)^{1-a}$. Consequently, we can write the $F_2$ series as follows:
\begin{eqnarray*}
\app2{a;1,1}{2,2}{x,y} \equal \sum_{i=0}^{\infty} \frac{(a)_i\,x^i}{(i+1)!}
\sum_{i=0}^{\infty} \frac{(a+i)_j\,y^j}{(j+1)!}\\
\equal \sum_{i=0}^{\infty} \frac{(a)_{i-1}\,x^i}{(i+1)!}\,\frac{(1-y)^{1-a-i}-1}y \\
\equal \frac{(1-y)^{1-a}}{(a-1)\,y}\sum_{i=0}^{\infty} \frac{(a-1)_{i}}{(i+1)!}\,\frac{x^i}{(1-y)^i} 
-\frac1{(a-1)\,y} \sum_{i=0}^{\infty} \frac{(a-1)_{i}\,x^i}{(i+1)!} \\
\equal \frac{(1-y)^{2-a}}{(a-1)\,(2-a)\,x\,y}
\left(1-\left(1-\frac{x}{1-y}\right)^{2-a} \right)-\frac{1-(1-x)^{2-a}}{(a-1)\,(2-a)\,x\,y}.
\end{eqnarray*}
The claimed formula follows. \qed\\

\noindent
Formula (\ref{eq:specf2}) is easy to specialize, if $\mbox{Re}\;a<2$:
\begin{eqnarray}
\app2{a;1,1\,}{2,2}{x,1}\equal \frac{1-(1-x)^{2-a}+(-x)^{2-a}}{(1-a)(2-a)\,x},\\
\app2{a;1,1\,}{2,2}{x,1-x}\equal \frac{1-(1-x)^{2-a}-x^{2-a}}{(1-a)(2-a)\,x\,(1-x)}.
\end{eqnarray}
Application of wrong formulas (\ref{we:ap2y1})--(\ref{we:ap2xy1}) evaluates the same two functions
to, respectively,
\begin{eqnarray}
\frac1{1-a}\,\hpg32{a,\,1,\,a-1}{2,\,a}{\,x} \equal
\frac1{1-a}\,\hpg21{1,\,a-1}{2}{\,x} \nonumber \\ \equal \frac{1-(1-x)^{2-a}}{(1-a)(2-a)\,x},\\
\frac{(1-x)^{-a}}{1-a}\,\hpg21{1,\,a-1}{2}{\frac{x}{x-1}}\equal \frac{1-(1-x)^{2-a}}{(1-a)(2-a)\,x\,(1-x)}.
\end{eqnarray}
Obviously, the numerator terms $(-x)^{2-a}$ or $-x^{2-a}$ are missing. 
The local exponent $2-a$ at $x=0$ for the respective third order Fuchsian equations does come into play. 

For the sake of completeness, here are expressions for the $F_2$ function in (\ref{eq:specf2})
with $a=1$ or $a=2$, from \cite{murleysaad}:
\begin{eqnarray}
\app2{1;1,1\,}{2,2}{x,y}\equal \frac{(1\!-\!x\!-\!y)\ln(1\!-\!x\!-\!y)-(1\!-\!x)\ln(1\!-\!x)-(1\!-\!y)\ln(1\!-\!y)}{x\,y},\\
\app2{2;1,1\,}{2,2}{x,y}\equal \frac{\ln(1-x)+\ln(1-y)-\ln(1-x-y)}{x\,y}.
\end{eqnarray}
Lemma \ref{eq:specf2} and these two formulas can be proved using 
the Euler type integral representation  \cite[9.4.(5)]{Srivastava85} for Appell's $F_2$ function.

\section{Univariate specializations of Appell's functions}
\label{sec:appell}

The paper \cite{AppelGhpg} proved that the following function pairs
satisfy the same ordinary differential  equations (of order $2$ or $3$):
\begin{enumerate}
\item $\displaystyle\app2{a;\;b_1,b_2}{c_1,c_2}{x,\,0}$ and $\displaystyle\hpg21{a,\,b_1}{c_1}{x}$;
\item $\displaystyle\app2{a;\;b_1,b_2}{c_1,c_2}{x,\,1}$ and 
$\displaystyle\hpg32{a,\,b_1,\,a-c_2+1}{c_1,\,a+b_2-c_2+1}{x}$;\\
\item $\displaystyle\app2{a;\;b_1,b_2}{c_1,c_2}{x,1-x}$ and
$\displaystyle(1-x)^{-a}\,\hpg32{a,\,c_1-b_1,\,a-c_2+1}{c_1,\,a+b_2-c_2+1}{\frac{x}{x-1}}$;\\
\item $\displaystyle\app3{a_1,a_2;\,b_1,b_2}{c}{x,\,0}$ and $\displaystyle\hpg21{a_1,\,b_1}{c}{x}$;
\item $\displaystyle\app3{a_1,a_2;\;b_1,b_2}{c}{x,\,1}$ and
$\displaystyle\hpg32{a_1,b_1,c-a_2-b_2}{c-a_2,c-b_2}{x}$;\\
\item $\displaystyle\app3{a_1,a_2;\;b_1,b_2}{c}{x,\frac{x}{x-1}}$ and\\
$\displaystyle x^{1-c}\,(1-x)^{a_2}\,\hpg32{1+a_1+a_2-c,1+b_1+a_2-c,1-b_2}
{1+a_1+a_2+b_1-c,1+a_2-b_2}{1-x}$;\\
\item $\displaystyle\app4{a;\;b}{c,a+b-c+\frac32}{x^2,(1-x)^2}$ 
and $\displaystyle\hpg21{2a,\,2b}{2c-1}{x}$;\\
\item $\displaystyle\app4{a;\,b}{c+\frac12,\frac12}{x^2,(1-x)^2}$ and 
$\displaystyle\hpg32{2a,\,2b,\,c}{a+b+\frac12,2c}{x}$;
\item $\displaystyle\app4{a;\;b}{c,a+b-c+1}{x^2,(1-x)^2}$
and $\displaystyle\hpg21{a,b\,}{c}{\,x}^2$;
\item $\displaystyle\app4{2c-\frac12;\,3c-1}{c+\frac12,c+\frac12}{x^2,(1-x)^2}$ and
$\displaystyle\hpg21{c,\,3c-1}{2c}{x}^2$;
\item $\displaystyle\app1{a;\;b_1,b_2}{c}{x,\,0}$ and $\displaystyle\hpg21{a,\,b_1}{c}{x}$;
\item $\displaystyle\app1{a;\;b_1,b_2}{c}{x,\,1}$ and
$\displaystyle\hpg21{a,\,b_1}{c-b_2}{x}$;
\item $\displaystyle\app1{a;\;b_1,b_2}{c}{x,\,x}$ and $\displaystyle\hpg21{a,\,b_1+b_2}{c}{x}$.
\end{enumerate}
The functions in {\em (i), (iv), (xi)} can be related straightforwardly as identical.
The equality sign can be put between the functions in {\em (xiii)} as well,  
as their power series at $x=0$ are comparable directly. 
The functions in {\em (v)} and {\em (xii)} are related in (\ref{eq:f31}) and (\ref{eq:f11}), respectively. 
The functions in {\em (vi)} have different behavior at the three singular points $x=0$, $x=1$, $x=\infty$, 
so they are not related by a two-term formula in general.
Naive identification of the $F_2$ functions in {\em (ii)} and {\em (iii)} is the topic of our main example.
The $F_4$ functions in {\em (vii)--(x)} cannot be directly identified by similar reasons.

An illustrative example is the function pair in {\em (ix)}. We can specialize Bailey's identity
\cite{Bailey33} 
\begin{equation} \label{eq:bailey}
\app4{a;\;b}{c,a\!+\!b\!-\!c\!+\!1}{x(1-y),y(1-x)}
=\hpg21{a,b}{c}{x} \hpg21{a,\;b}{a\!+\!b\!-\!c\!+\!1}{y}
\end{equation}
to 
\begin{equation} \label{eq:bailey2}
\app4{a;\;b}{c,a\!+\!b\!-\!c\!+\!1}{x^2,(1-x)^2}
=\hpg21{a,b}{c}{x} \hpg21{a,\;b}{a\!+\!b\!-\!c\!+\!1}{1-x}.
\end{equation}
Both sides of this equality and satisfy the symmetric tensor square Fuchsian equation
for $\hpg21{\!a,\,b}c{z}^2$. For $x\in[0,1]$, the specialization is on the border of the convergence 
region of the $F_4(u,v)$ function.  To get a local expression at $x=0$, one can use
a connection formula in \cite[Theorem 2.3.2]{specfaar} to transform the $\hpgo21(1-x)$ function, 
if $\mbox{Re }c<1$ and $\mbox{Re}(c-a-b)>0$:
\begin{eqnarray} \label{eq:bailey2a}
\app4{a;\;b}{c,a\!+\!b\!-\!c\!+\!1}{x^2,(1-x)^2} =
\frac{\Gamma(a+b-c+1)\Gamma(1-c)}{\Gamma(a-c+1)\Gamma(b-c+1)}
\,\hpg21{a,b}{c}{x}^2 \qquad \nonumber \\
+\frac{\Gamma(a\!+\!b\!-\!c\!+\!1)\Gamma(c\!-\!1)}{\Gamma(a)\Gamma(b)}
\,x^{1-c}\,\hpg21{a,b}{c}{x}\hpg21{a-c+1,b-c+1}{2-c}{x}.
\end{eqnarray}

Consider next the function pair in {\em (vii)}. It satisfies a second order Fuchsian equation, 
and we may relate its various hypergeometric solutions by evaluating at two points, 
say $x=0$ and $x=1$. In this way we get
\begin{eqnarray} \label{eq:f4sa}
\app4{a;\;b}{c,a+b-c+\frac32}{x^2,(1-x)^2} = \hspace{-150pt} \nonumber \\
&&\frac{\Gamma(a+b-c+\frac32)\Gamma(\frac32-c)}{\Gamma(a-c+\frac32)\Gamma(b-c+\frac32)}\,
\frac{\cos\pi c\,\cos\pi(c-a-b)}{\cos\pi a\,\cos\pi b}\,\hpg21{2a,\,2b}{2c-1}{x}+\nonumber \\
&& \frac{\Gamma(c)\Gamma(c-a-b)}{\Gamma(c-a)\Gamma(c-b)}\,
\frac{\sin\pi c\,\sin\pi(c-a-b)}{\cos\pi a\,\cos\pi b}\,
\hpg21{2a,\,2b}{2a\!+\!2b\!-\!2c\!+\!2}{1-x}
\end{eqnarray}
on the interval $[0,1]$, if $\mbox{Re }c<1$ and $\mbox{Re}(c-a-b)>\frac12$.
We can use a connection formula in \cite[Theorem 2.3.2]{specfaar} 
and eliminate the $\hpgo21(1-x)$ function here as well:
\begin{eqnarray} \label{eq:f4sb}
\!\app4{a;\;b}{\!c,a\!+\!b\!-\!c\!+\!\frac32}{x^2,(1-x)^2}\! \equal 
\frac{\Gamma\left(a+b-c+\frac32\right)\Gamma\left(\frac32-c\right)}
{\Gamma\left(a-c+\frac32\right)\Gamma\left(b-c+\frac32\right)}\,
\hpg21{2a,\,2b}{2c-1}{x} +\nonumber \\
&&\!\frac{\Gamma(c)\Gamma\!\left(a\!+\!b\!-\!c\!+\!\frac32\right)\!\Gamma(2c\!-\!2)}
{\Gamma(c\!-\!a)\,\Gamma(c\!-\!b)\,\Gamma(2a)\,\Gamma(2b)}\,
\frac{2^{2a+2b-2c+1}\sqrt{\pi}\sin\pi c}{\cos\pi a\,\cos\pi b}\quad\nonumber \\
&&\!\times\, x^{2-2c}\,\hpg21{2a-2c+2,\,2b-2c+2}{3-2c}{\,x}.
\end{eqnarray}
Similar identities for the cases {\em (viii), (x)} may involve in total four different solutions
of the same Fuchsian equation; the coefficients can be derived if evaluation at 3 points is known,
say $x=0$, $x=1$, $x=\frac12$. Verification of (\ref{eq:f4sa})--(\ref{eq:f4sb}) at $x=0$ and $x=1$ 
is based on the following trigonometric identities:
\begin{eqnarray}
\cos\pi a\,\cos\pi b\equal \cos\pi(c-a)\,\cos\pi(c-b)+\sin\pi c\,\sin\pi(c-a-b)\\
\equal \sin\pi(c-a)\,\sin\pi(c-b)+\cos\pi c\,\cos\pi(c-a-b),\\
\sin\pi a\,\sin\pi b\equal \sin\pi(c-a)\,\sin\pi(c-b)-\sin\pi c\,\sin\pi(c-a-b)\\
\equal \cos\pi(c-a)\,\cos\pi(c-b)-\cos\pi c\,\cos\pi(c-a-b).
\end{eqnarray}

All Appell's functions in {\em (i)--(xiii)} are specialized to their singular (aka branching) curves. 
But the same routines apply to univariate specializations of Appell's functions to outside the singular locus. For example, \cite{AppelGhpg} proves that the following function pairs satisfy the same
Fuchsian equations of order 2:
\begin{enumerate}
\item[{\em (xiv)}] $\displaystyle\app1{a;\;2b, a-b}{1+b}{x,\,x^2}$ and
$\displaystyle(1-x)^{-2a}\,\hpg21{a,\;\frac12}{1+b}{-\frac{4x}{(x-1)^2}}$;
\item[{\em (xv)}]  $\displaystyle\app2{a;\;b_1,b_2}{2b_1,2b_2}{x,2-x}$ and
$\displaystyle(x-2)^{-a}\,\hpg21{\frac{a}2,\frac{a+1}2-b_2}{b_1+\frac12}{\frac{x^2}{(2-x)^2}}$.
\end{enumerate}
The functions in {\em (xiv}) can be identified as equal, since their series at $x=0$ identify straightforwardly. The specialization in {\em (xv)} is onto a line way outside the region of convergence
of the $F_2(u,v)$ series. To relate the $F_2(x,2-x)$ and $\hpgo21$ solutions, one has to consider analytic continuation of the $F_2(u,v)$ function scrupulously, or assume the $F_2(u,v)$ series terminating in at least one variable.

\section{Claussen-type example}

Let us consider the following bivariate series:
\begin{eqnarray} \label{eq:gap211}
F^{2:1;1}_{1:1;1}\!\left(\left. {a; b; p_1, p_2\atop c; q_1, q_2} \right| x,y\right) \equal
\sum_{i=0}^{\infty} \sum_{j=0}^{\infty}
\frac{(a)_{i+j}\,(b)_{i+j}\,(p_1)_i\,(p_2)_j}{(c)_{i+j}\,(q_1)_i\,(q_2)_j\;i!\,j!}\,x^i\,y^j.
\end{eqnarray}
It is a special case of {\em Kamp\'e de F\'eriet} series.
If $a=c$, it reduces to Appell's $F_2$ series. In \cite{gclausen}, it is proved that
\begin{equation} \label{eq:gap211a}
\hpg21{a,\,b}{c}{\,z}^2 \quad\mbox{and}\quad
F^{2:1;1}_{1:1;1}\!\left(\left. {2a;\,2b;\;c-\frac12,\,a+b-c+\frac12 
\atop a+b+\frac12;\, 2c-1, 2a+2b-2c+1} \right|  z,1-z\right)
\end{equation}
satisfy the same Fuchsian ordinary differential equation of order 3.
This fact generalizes Clausen's identity \cite{clausen28}:
\begin{equation} \label{eq:clausen}
\hpg21{a,\;b}{a+b+\frac12}{\,z}^2=\hpg32{2a,\,2b,\,a+b}{2a+2b,a+b+\frac12}{\,z}.
\end{equation}

The $F^{2:1;1}_{1:1;1}$ double series cannot be directly re-expanded as 
Taylor-series at neither $z=0$ nor $z=1$, for the similar reasons as (\ref{eq:app2x0}). 
Therefore the functions in (\ref{eq:gap211a}) are generally not related by a two-term identity.
In particular, the evaluations at $z=0$ and $z=1$ are not consistent. We have
\begin{eqnarray} \label{eq:kdf1}
F^{2:1;1}_{1:1;1}\!\left(\left. {2a;\,2b;\;c-\frac12,\,a+b-c+\frac12 
\atop a+b+\frac12;\, 2c-1, 2a+2b-2c+1} \right|  0,1\right)\equal \nonumber\\
&&\hspace{-110pt}
\frac{\Gamma\left(\frac12\right)\Gamma\left(a+b+\frac12\right)\Gamma(1-c)\,\Gamma(1+a+b-c)}
{\Gamma\left(a+\frac12\right)\Gamma\left(b+\frac12\right)\Gamma(1+a-c)\,\Gamma(1+b-c)},\\
\label{eq:kdf2} F^{2:1;1}_{1:1;1}\!\left(\left. {2a;\,2b;\;c-\frac12,\,a+b-c+\frac12 
\atop a+b+\frac12;\, 2c-1, 2a+2b-2c+1} \right|  1,0\right)\equal \nonumber\\
&&\hspace{-110pt}
\frac{\Gamma\left(\frac12\right)\Gamma\left(a+b+\frac12\right)\Gamma(c)\,\Gamma(c-a-b)}
{\Gamma\left(a+\frac12\right)\Gamma\left(b+\frac12\right)\Gamma(c-a)\,\Gamma(c-b)},
\end{eqnarray}
by Watson's $\hpgo32(1)$ summation formula \cite[Theorem 3.5.5\em (i)]{specfaar}, 
if $\mbox{Re } c<1$ and \mbox{$\mbox{Re}(c-a-b)>1$}. 
But
\begin{eqnarray*}
\hpg21{a,\;b}{a+b+\frac12}{\,0}^2=1,\qquad
\hpg21{a,\;b}{a+b+\frac12}{\,1}^2=
\frac{\Gamma(c)^2\,\Gamma(c-a-b)^2}{\Gamma(c-a)^2\,\Gamma(c-b)^2}.
\end{eqnarray*}
The ration of these two squares of $\hpgo21$ functions is generally {\em not} equal to the ratio
of the right-hand sides of (\ref{eq:kdf1}) and (\ref{eq:kdf2}), so the two functions in (\ref{eq:gap211a})
cannot be proportional.

However, the naive evaluation still gives the non-branching series part of the 
$F^{2:1;1}_{1:1;1}$ function at $z=0$ right. The $F^{2:1;1}_{1:1;1}$ function must be
the following linear combination of local solutions at $z=0$ of the third order Fucshian equation:
\begin{eqnarray} \label{ge:clausen}
F^{2:1;1}_{1:1;1}\!\left(\left. {2a;\,2b;\;c-\frac12,\,a+b-c+\frac12 
\atop a+b+\frac12;\, 2c-1, 2a+2b-2c+1} \right|  z,1-z\right)= \hspace{-204pt} \nonumber \\ 
&& \frac{\Gamma\left(\frac12\right)\Gamma\left(a+b+\frac12\right)\Gamma(1-c)\,\Gamma(1+a+b-c)}
{\Gamma\left(a+\frac12\right)\Gamma\left(b+\frac12\right)\Gamma(1+a-c)\,\Gamma(1+b-c)}\,
\hpg21{a,\;b}{c}{\,z}^2 \nonumber\\ 
&& + \mbox{ the term with } z^{1-c}\,\hpg21{a,\;b}{c}{\,z}\hpg21{1+a-c,1+b-c}{2-c}{z} \nonumber\\
&& + \mbox{ the term with } z^{2-2c}\,\hpg21{1+a-c,1+b-c}{2-c}{z}^2.
\end{eqnarray}
This form can be used to prove the identity \cite[(10)]{gclausen}:
\begin{eqnarray}  \label{ge:clausen2}
\hpg21{a,\;b}{a+b+n+\frac12}{\,z}^2\! \equal
\frac{(\frac12)_n\,(a+b+\frac12)_n}{(a+\frac12)_n\,(b+\frac12)_n} \times \nonumber\\ 
& & F^{2:1;1}_{1:1;1}\!\left(\left. {2a;\,2b;\;a+b+n,\,-n 
\atop a+b+\frac12;\, 2a+2b+2n, -2n} \right|  z,1-z\right).
\end{eqnarray}
Here the series in the second argument $1-z$ is understood to be terminating with the power $(1-z)^n$,
so the series at $z=0$ of the right-hand side is well defined, and branching terms
are not present. Therefore the $\hpgo21(z)^2$ and $F^{2:1;1}_{1:1;1}(z,1-z)$ functions
must differ by a constant factor, which can be determined after evaluation at $z=1$
as in \cite{gclausen}. But we may start by taking limit $c\to a+b+n+\frac12$ in (\ref{ge:clausen}).
Since
\begin{equation}
\lim_{\varepsilon\to 0} \frac{(\varepsilon-n)_{2k+1}}{(2\varepsilon-2n)_{2n+1}}=
\frac{(-n)_n\,n!}{(-2n)_{2n}\cdot2}=\frac{(-1)^n\,(n!)^2}{2\cdot(2n)!}
=\frac{(-1)^n\,n!}{2^{2n+1}(\frac12)_n},
\end{equation}
 we have
\begin{eqnarray} \label{ge:altern1}
&&\hspace{-50pt} \frac{\Gamma\left(\frac12\right)\Gamma\left(a+b+\frac12\right)
\Gamma(\frac12-a-b-n)\,\Gamma(\frac12-n)}
{\Gamma\left(a+\frac12\right)\Gamma\left(b+\frac12\right)\Gamma(\frac12-a-n)\,\Gamma(\frac12-b-n)}
\,\hpg21{a,\;b}{a+b+n+\frac12}{\,z}^2 \hspace{-272pt} \nonumber\\ 
&& \hspace{-50pt}+ \mbox{ generally branching power series terms} \hspace{-272pt} \nonumber\\ 
\equal F^{2:1;1}_{1:1;1}\!\left(\left. {2a;\,2b;\;a+b+n,\,-n 
\atop a+b+\frac12;\, 2a+2b+2n, -2n} \right|  z,1-z\right)+ \nonumber\\
&&\frac{(-1)^n\,n!}{2^{2n+1}\,(\frac12)_n}
\frac{(2a)_{2n+1}(2b)_{2n+1}}{\,(2n+1)!\,(a+b+\frac12)_{2n+1}}\,(1-z)^{2n+1}\times \nonumber\\
&&F^{2:1;1}_{1:1;1}\!\left(\left. {2a+2n+1;\,2b+2n+1;\;a+b+n,\,n+1 
\atop a+b+2n+\frac32;\, 2a+2b+2n, 2n+2} \right|  z,1-z\right).
\end{eqnarray}
The latter $F^{2:1;1}_{1:1;1}$ function can be written as 
\begin{eqnarray*}
&&\frac{\Gamma\left(\frac12\right)\Gamma\left(a+b+2n+\frac32\right)
\Gamma(\frac12-a-b-n)\Gamma(\frac32+n)}
{\Gamma\left(a+n+1\right)\Gamma\left(b+n+1\right)\Gamma(1-a)\,\Gamma(1-b)}\,
\hpg21{a+n+\frac12,b+n+\frac12}{a+b+n+\frac12}{\,z}^2 \nonumber\\
&&+ \mbox{ generally branching power series terms}.
\end{eqnarray*}
Besides,
\begin{eqnarray*}
(1-z)^{n+1}\,\hpg21{a+n+\frac12,b+n+\frac12}{a+b+n+\frac12}{\,z}^2=
\hpg21{a,\;b}{a+b+n+\frac12}{\,z}^2
\end{eqnarray*}
by Euler's transformation \cite[(2.2.7)]{specfaar}. We also use $(2n+1)!=2^{2n+1}n!\,(\frac12)_{n+1}$ 
and Euler's reflection formula \cite[(1.2.1)]{specfaar} for the $\Gamma$-function
to rewrite identity (\ref{ge:altern1}) as follows:
\begin{eqnarray} 
&&\hspace{-60pt}\frac{\cos\pi a\,\cos\pi b}{\cos\pi(a+b)}\,
\frac{(\frac12-a-n)_n\,(\frac12-b-n)_n}{(\frac12-n)_n\,(\frac12-a-b-n)_n}\,
\hpg21{a,\;b}{a+b+n+\frac12}{\,z}^2 \hspace{-272pt}  \nonumber\\ 
&&\hspace{-60pt}+ \mbox{ generally branching power series terms} \nonumber\\ 
\equal F^{2:1;1}_{1:1;1}\!\left(\left. {2a;\,2b;\;a+b+n,\,-n 
\atop a+b+\frac12;\, 2a+2b+2n, -2n} \right|  z,1-z\right)+ \nonumber\\
&&\frac{(-1)^n}{2^{4n+2}\,(\frac12)_n\,(\frac12)_{n+1}}
\frac{(2a)_{2n+1}(2b)_{2n+1}}{(a+b+\frac12)_{2n+1}}\,\frac{\sin\pi a\,\sin\pi b}{\cos\pi(a+b)} 
\times \nonumber\\ &&
\frac{(\frac12)_{n+1}\,(a+b+\frac12)_{2n+1}}{(a)_{n+1}\,(b)_{n+1}\,(\frac12-a-b-n)_n}
\hpg21{a,\;b}{a+b+n+\frac12}{\,z}^2 \nonumber \\
&&+ \mbox{ generally branching power series terms}.
\end{eqnarray}
The generally branching terms must cancel.
We collect the coefficients to the $\hpgo21(z)^2$ function on the left-hand side.
The resulting factor is
\begin{eqnarray*}
\frac{\cos\pi a\,\cos\pi b}{\cos\pi(a+b)}\,
\frac{(a+\frac12)_n\,(b+\frac12)_n}{(\frac12)_n\,(a+b+\frac12)_n}-
\frac{\sin\pi a\,\sin\pi b}{\cos\pi(a+b)}
\frac{(2a)_{2n+1}\,(2b)_{2n+1}}{2^{4n+2}(\frac12)_n(a+b+\frac12)_n(a)_{n+1}(b)_{n+1}}\\
=\frac{(a+\frac12)_n\,(b+\frac12)_n}{(\frac12)_n\,(a+b+\frac12)_n},
\end{eqnarray*}
and formula (\ref{ge:clausen2}) is proved.

\bibliographystyle{alpha}
\bibliography{../Transformations/hypergeometric}

\begin{thebibliography}{AAR99}

\bibitem[AAR99]{specfaar}
G.E. Andrews, R.~Askey, and R.~Roy.
\newblock {\em Special Functions}.
\newblock Cambridge Univ. Press, Cambridge, 1999.

\bibitem[Bai33]{Bailey33}
W.~N. Bailey.
\newblock A reducible case of the fourth type of appell's hypergeometric
  functions of two variables.
\newblock {\em Quart. J. Math. (Oxford)}, 4:305--308, 1933.

\bibitem[Cla28]{clausen28}
T.~Clausen.
\newblock Ueber die {F}\"alle, wenn die {R}eihe von der {F}orm
  $y=1+\frac{\alpha}{1}\!\cdot\!\frac{\beta}{\gamma}x+
  \frac{\alpha\cdot\alpha+1}{1\cdot2}\!\cdot\!
  \frac{\beta\cdot\beta+1}{\gamma\cdot\gamma+1}x^2+$etc. ein quadrat von der
  {F}orm $z=1+\frac{\alpha'}{1}\cdot\frac{\beta'}{\gamma'}\!\cdot\!
  \frac{\delta'}{\epsilon'}x+\frac{\alpha'\cdot\alpha'+1}{1\cdot2}
  \!\cdot\!\frac{\beta'\cdot\beta'+1}{\gamma'\cdot\gamma'+1}
  \!\!\cdot\frac{\delta'\cdot\delta'+1}{\epsilon'\cdot\epsilon'+1}x^2+$etc.
  hat.
\newblock {\em J. Reine Ang. Math.}, 3:89--91, 1828.

\bibitem[MS08]{murleysaad}
J.~Murley and N.~Saad.
\newblock Tables of the {A}ppell hypergeometric functions ${F}_2$.
\newblock Available at {\sf http://arxiv.org/abs/0809.5203}., 2008.

\bibitem[SK85]{Srivastava85}
H.~M. Srivastava and P.~W. Karlsson.
\newblock {\em Multiple Gaussian Hypergeometric Series}.
\newblock Ellis Horwood Ltd., 1985.

\bibitem[Vid09a]{gclausen}
R.~Vid\=unas.
\newblock A generalization of {C}lausen's identity.
\newblock Available at {\sf http://arxiv.org/abs/0906.1862}, 2009.

\bibitem[Vid09b]{AppelGhpg}
R.~Vid\=unas.
\newblock Specialization of {A}ppell's functions to univariate hypergeometric
  functions.
\newblock {\em Journal of Mathematical Analysis and Applications},
  355:145--163, 2009.
\newblock Available at {\sf http://arxiv.org/abs/0804.0655}.

\end{thebibliography}

\end{document}